\documentclass[12pt]{amsart}
\usepackage{graphicx}
\usepackage{times}
\usepackage{mathptm}

\begin{document}
\newtheorem{theorem}{Theorem}[section]
\newtheorem{prop}[theorem]{Proposition}
\newtheorem{lemma}[theorem]{Lemma}
\newtheorem{claim}[theorem]{Claim}
\newtheorem{num}{Counting Rule}
\newtheorem{cor}[theorem]{Corollary}
\newtheorem{defin}[theorem]{Definition}
\newtheorem{example}[theorem]{Example}
\newtheorem{xca}[theorem]{Exercise}
\newcommand{\map}{\mbox{$\rightarrow$}}
\newcommand{\aaa}{\mbox{$\alpha$}}
\newcommand{\Aaa}{\mbox{$\mathcal A$}}
\newcommand{\bbb}{\mbox{$\beta$}}
\newcommand{\ccc}{\mbox{$\mathcal C$}}
\newcommand{\ddd}{\mbox{$\delta$}} 
\newcommand{\Ddd}{\mbox{$\Delta$}}
\newcommand{\Fff}{\mbox{$\mathcal F$}}  
\newcommand{\Ggg}{\mbox{$\Gamma$}}
\newcommand{\ggg}{\mbox{$\gamma$}}
\newcommand{\kkk}{\mbox{$\kappa$}}
\newcommand{\lll}{\mbox{$\lambda$}}
\newcommand{\Lll}{\mbox{$\Lambda$}}
\newcommand{\mlp}{\mbox{$\mu^{+}_{l}$}}
\newcommand{\ml}{\mbox{$\mu_{l}$}}
\newcommand{\mr}{\mbox{$\mu_{r}$}}
\newcommand{\mlpm}{\mbox{$\mu_{l}^{\pm}$}}
\newcommand{\mrpm}{\mbox{$\mu_{r}^{\pm}$}}
\newcommand{\mlm}{\mbox{$\mu_{l}^{-}$}}
\newcommand{\mrp}{\mbox{$\mu_{r}^{+}$}}
\newcommand{\mrm}{\mbox{$\mu_{r}^{-}$}}
\newcommand{\mm}{\mbox{$\mu^-$}}
\newcommand{\mpm}{\mbox{$\mu^{\pm}$}}
\newcommand{\mpp}{\mbox{$\mu^+$}}
\newcommand{\mt}{\mbox{$\mu^{t}$}}
\newcommand{\mb}{\mbox{$\mu_{b}$}}
\newcommand{\mz}{\mbox{$\mu^{\bot}$}}
\newcommand{\mpq}{\mbox{$\mu^{(p,q)}$}}
\newcommand{\omp}{\mbox{$0_{-}^{+}$}}
\newcommand{\oa}{\mbox{$\overline{a}$}}
\newcommand{\ob}{\mbox{$\overline{b}$}}
\newcommand{\opm}{\mbox{$0_{+}^{-}$}}
\newcommand{\opp}{\mbox{$0_{+}^{+}$}}
\newcommand{\Pt}{\mbox{$\tilde{P}$}}
\newcommand{\rrr}{\mbox{$\rho$}} 
\newcommand{\rz}{\mbox{$\rho^{\bot}$}}
\newcommand{\rp}{\mbox{$\rho^+$}}
\newcommand{\rmm}{\mbox{$\rho^-$}}
\newcommand{\rpq}{\mbox{$\rho^{(p,q)}$}}
\newcommand{\sss}{\mbox{$\sigma$}} 
\newcommand{\Sss}{\mbox{$\mathcal S$}} 
\newcommand{\sm}{\mbox{$\sigma^-$}}
\newcommand{\sz}{\mbox{$\sigma^{\bot}$}} 
\newcommand{\spm}{\mbox{$\sigma^{\pm}$}}
\newcommand{\spp}{\mbox{$\sigma^+$}}
\newcommand{\st}{\mbox{$\sigma^{t}$}}
\newcommand{\Ss}{\mbox{$\Sigma$}}
\newcommand{\Th}{\mbox{$\Theta$}} 
\newcommand{\ttt}{\mbox{$\tau$}} 
\newcommand{\bdd}{\mbox{$\partial$}}
\newcommand{\zzz}{\mbox{$\zeta$}}
\newcommand{\qb} {\mbox{$Q_B$}}
\newcommand{\inter}{\mbox{${\rm int}$}}

\title[] {Thinning genus two Heegaard spines in $S^{3}$}

\author{Martin Scharlemann}
\address{\hskip-\parindent
        Mathematics Department\\
        University of California\\
        Santa Barbara, CA 93106\\
        USA}
\email{mgscharl@math.ucsb.edu}

\author{Abigail Thompson}
\address{\hskip-\parindent
        Mathematics Department\\
        University of California\\
        Davis, CA 95616\\
        USA}
\email{thompson@math.ucdavis.edu}

\date{\today} 
\thanks{Research supported in part by NSF grants and by the generous 
support of RIMS, Kyoto}

\begin{abstract} We study trivalent graphs in $S^{3}$ whose closed 
complement is 
a genus two handlebody.  We show that such a graph, when put in thin 
position, has a level edge connecting two vertices.

\end{abstract}
\maketitle
  
\section{Introduction}

We briefly review the terminology of Heegaard splittings, referring
the reader to \cite{Sc} for a more complete description.  A {\em 
Heegaard splitting} of a closed $3$-manifold $M$ is a division of $M$ 
into two handlebodies by a connected closed surface, called the 
Heegaard surface or the splitting surface.  A {\em spine} for a 
handlebody $H$ is a graph $\Ggg \subset interior(H)$ so that $H$ is a 
regular neighborhood of $\Ggg$.  A {\em Heegaard spine} in $M$ is a 
graph $\Ggg \subset M$ whose regular neighborhood $\eta(\Ggg)$ has 
closed complement a handlebody.  Equivalently, $\bdd\eta(\Ggg)$ is a 
Heegaard surface for $M$.  We say that $\Ggg$ is of genus $g$ if 
$\bdd\eta(\Ggg)$ is a surface of genus $g$.

Any two spines of the same handlebody are equivalent under edge 
slides (see \cite{ST1}).  It's a theorem of Waldhausen \cite{Wa} (see 
also \cite{ST2}) that any Heegaard splitting of $S^{3}$ can be 
isotoped to a standard one of the same genus.  An equivalent 
statement, then, is that any Heegaard spine for $S^{3}$ can be made 
planar by a series of edge slides.

On the other hand, without edge slides, Heegaard spines in $S^{3}$ can
be quite complicated, even for genus as low as two.  For example, let
$L$ be a $2$-bridge knot or link in bridge position and $\ggg$ be a
level arc that connects the top two bridges.  Then it's easy to see
that the graph $L \cup \ggg$ is a Heegaard spine since, once $\ggg$ is
attached, the arcs of $\Ggg$ descending from $\ggg$ can be slid around
on $\ggg$ until the whole graph is planar.  More generally, a knot or
link is called {\em tunnel number one} if the addition of a single arc
turns it into a Heegaard spine.  For Heegaard spines constructed in
this way, it was shown in \cite{GST} that the picture for the
two-bridge knot is in some sense the standard picture.  That is, if
$L$ is a tunnel number one knot or link put in minimal bridge
position, and $\ggg$ is an unknotting tunnel, then $\ggg$ may be slid
on $L$ and isotoped in $S^{3}$ until it is a level arc.  The ends of
$\ggg$ may then be at the same point of $K$ (so $\ggg$ becomes an
unknotted loop) or at different points (so $\ggg$ becomes a level
edge).  It can even be arranged that, when $\ggg$ is level, the ends
of $\ggg$ lie on one or two maxima (or minima).  Finally, in
\cite{GST} the notion of {\em width} for knots was extended to
trivalent graphs, and it was shown that this picture of $L \cup \ggg$
is in some sense natural with respect to this measure of complexity. 
Specifically, if $\ggg$ is slid and isotoped to make the graph $\Ggg$
as thin as possible without moving $L$, then $\ggg$ will be made
level.

This raises a natural question.  As we've seen, choosing $\ggg$ to 
make $L \cup \ggg$ as thin as possible reveals that $\ggg$ can 
actually be made level.  So suppose $\Ggg$ is an arbitrary Heegaard 
spine of $S^{3}$ (but trivalent so the notion of thin position is 
defined) and we allow no edge slides at all.  Suppose a height 
function is given on $S^{3}$ and we isotope $\Ggg$ in $S^{3}$ to make 
it as thin as possible.  What can then be said about the positioning 
of $\Ggg$?  We will answer the question for genus two Heegaard spines 
by showing this: once a trivalent genus $2$ Heegaard spine $\Ggg$ is 
put in thin position, some simple edge (that is, an edge not a loop) 
will be level. It is an intriguing question whether there is any 
analogous result for higher genus Heegaard spines.

Here is an outline of the rest of the paper.  In Section 2, we give
some definitions and we prove prove a preliminary proposition
(Proposition \ref{prop:mori}) generalizing a theorem of Morimoto, and
a preliminary lemma (Lemma \ref{lemma:bridgeye}) for eyeglass graphs. 
In section 3 we state and prove the two main theorems of the paper
(Theorems \ref{theorem:bridgeye} and \ref{theorem:thineye}) together
with Corollary \ref{cor:leveledge} which gives the result stated in
the abstract.  In section 4 we state and prove a technical lemma
(Lemma \ref{lemma:tech}) needed in the proof of Theorem
\ref{theorem:thineye}.

\section{Preliminaries}

\begin{figure}
\centering
\includegraphics [width=0.65\textwidth]{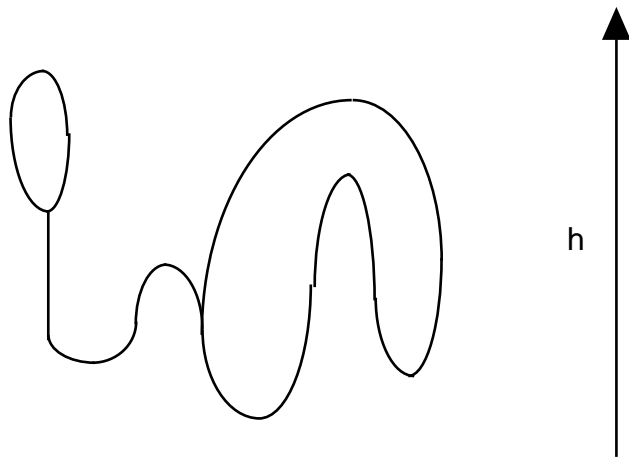}
\caption{extended bridge position}  \label{extbridge}
\end{figure}

\begin{defin} Let $\Ggg \subset S^{3}$ be a trivalent graph.  Suppose
a height function is defined on $S^{3}$.  A cycle in $\Ggg$ is {\em
vertical} if it has exactly one minimum and one maximum.  $\Ggg$ is in
{\em bridge position} if every minimum lies below every maximum.  A
{\em regular} minimum or maximum is one that does not occur at a
vertex.  A trivalent graph is in {\em extended bridge position}
(Figure \ref{extbridge}) if any minimum lying above a regular maximum
(resp.  maximum lying below a regular minimum) is a $Y$-vertex at the
minimum (resp.  $\lll$-vertex at the maximum) of a vertical cycle.
\end{defin}

\begin{defin} Suppose $\Ggg \subset S^{3}$ is a trivalent graph and 
$H$ is a regular neighborhood.  Let $\mu_{1}$, $\mu_{2}$ be two 
meridians of $H$ corresponding to points $p_{1}, p_{2}$ on $\Ggg$.  
Then a path $\aaa$ between the $\mu_{i}$ is {\em regular} if it is 
parallel in $H - (\mu_{1} \cup \mu_{2})$ to an embedded path in 
$\Ggg$.  That is, it intersects each meridian of $H$ in at most one 
point.
\end{defin}

\begin{figure}
\centering
\includegraphics [width=0.65\textwidth]{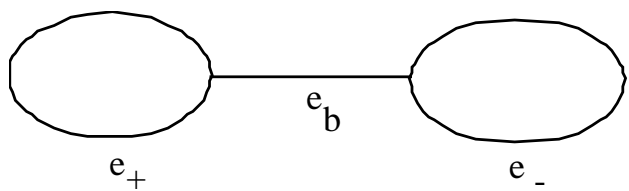}
\caption{eyeglass graph}  \label{eyeglass}
\end{figure}

\begin{defin} An {\em eyeglass graph} (Figure \ref{eyeglass}) is a
graph consisting of two cycles $e_{\pm}$ connected by an edge $e_{b}$,
called the bridge edge.
\end{defin}

We extend a theorem of Morimoto \cite{Mo} that extends earlier work 
of Gordon-Reid \cite{GR}:

\begin{prop}\label{prop:mori}  
Let $\Ggg$ be a trivalent Heegaard spine in $S^{3}$ whose regular
neighborhood $H$ is a genus two handlebody.  Suppose $Q$ is a
collection of spheres in general position with respect to $\Ggg$, so 
$Q$ intersects $H$ in a collection of meridians, each corresponding
to a point in $\Ggg \cap Q$.  Suppose $Q - H$ is incompressible in the
complement of $H$ and no component is a disk.  Then either:

\begin{enumerate}  

\item each component of $Q \cap \bdd H$ is a non-separating curve and
each component of $Q - H$ is parallel in $S^3 - H$ to a component of
$\bdd H - Q$

\item each component of $Q \cap \bdd H$ is a separating curve, and
each component of $Q - H$ is parallel in $S^3 - H$ to a component of
$\bdd H - Q$.  (Each component of $Q - H$ is then necessarily an
annulus).

\item $Q \cap \bdd H$ contains both separating and non-separating
curves.  Then there is a disk $F \subset S^{3}$ whose interior is
disjoint from $H \cup Q$ and $\bdd F = \aaa \cup \bbb$, where $\aaa
\subset \bdd H$, $\bbb \subset Q$.  Either

\begin{enumerate}
\item  $\aaa$ is a regular path on $\bdd H$ that is disjoint from 
some meridian corresponding to a point in $e_{b}$ or

\item $\aaa$ has both ends at the same separating meridian and 
intersects some non-separating meridian in exactly one point.

\end{enumerate}
\end{enumerate}
\end{prop}

{\bf Remark:} Of course, unless $\Ggg$ is an eyeglass whose bridge
edge is intersected by $Q$, only the first possibility is relevant. 
Notice also that in case 1) or case 2) then automatically 
there is a disk as described in item (3a).

\begin{proof} 
The first two cases are proven by Morimoto \cite{Mo}.  
So we assume $\Ggg$ is an eyeglass graph.  The proof will be by 
induction on $Q \cap e_{b}$; when $Q \cap e_{b} = \emptyset$ the 
result follows from case 1), so we assume $Q \cap e_{b} \neq \emptyset$.  

Let $E$ be an essential disk in the closed complement of $H$.  We can
assume that some component of $Q \cap H$ is a separating meridian, or
else item 1) would apply.  We can assume that $E \cap Q \neq
\emptyset$ or else some component of $Q$ with a separating meridian
would be compressible.  Let $E_{0}$ be an outermost arc of $E$ cut off
by $Q$.  Let $\aaa=E_{0} \cap \bdd H$, $\bbb=E_{0} \cap Q$.  We can
assume there are no disks of intersection between $E_{0}$ and $Q$
since $Q$ is incompressible.  If $\aaa$ connects distinct meridians of
$H$ we are done, for $\aaa$ is disjoint from the meridian
corresponding to any point in $Q \cap e_{b}$, so $E_{0}$ serves for
$F$ in (3a).  So we will suppose both ends of $\aaa$ lie at the same
meridian $\mu$ of $Q$.  A counting argument on the number of
intersection points between $\bdd E$ and the three natural meridian
curves on $\bdd \Ggg$ shows that $\mu$ cannot be non-separating.

So suppose $\mu$ is separating and $\aaa$ intersects a meridian of
$e_{+}$ non-trivially.  Join the ends of $\aaa$ together on $\mu$ to
get a closed curve $\aaa_{+}$ lying on the boundary of a solid torus
(essentially a neighborhood of $e_{+}$) and bounding a disk in its
complement (the disk is the union of $E_{0}$ and a disk in $Q$). 
Hence $\aaa_{+}$ is a longitude and we have item 3b); 
a meridian of $e_{+}$ is the meridian intersected once.

The interesting case is when $\mu$ is separating and $\aaa$ is a
``wave'' at an end of $e_b$, that is, $\aaa$ is disjoint from a
meridian of each cycle (Figure \ref{fig:wave}).  In this case, modify $Q$
by ``splitting'' the end of $e_{b}$ to which $\bdd E_{0}$ is incident. 
Equivalently, push out that meridian of $Q$ past the end of $e_{b}$ so
that it splits into two meridians of, say, $e_{-}$ (Figure \ref{Q¹}). 
Call the new collection of spheres $Q'$.  The splitting converts
$E_{0}$ into a compressing disk for $Q'$.  Let $Q_{0}$ be
the collection of spheres obtained by compressing $Q'$ along that
disk.

\begin{figure}
\centering
\includegraphics [width=0.65\textwidth]{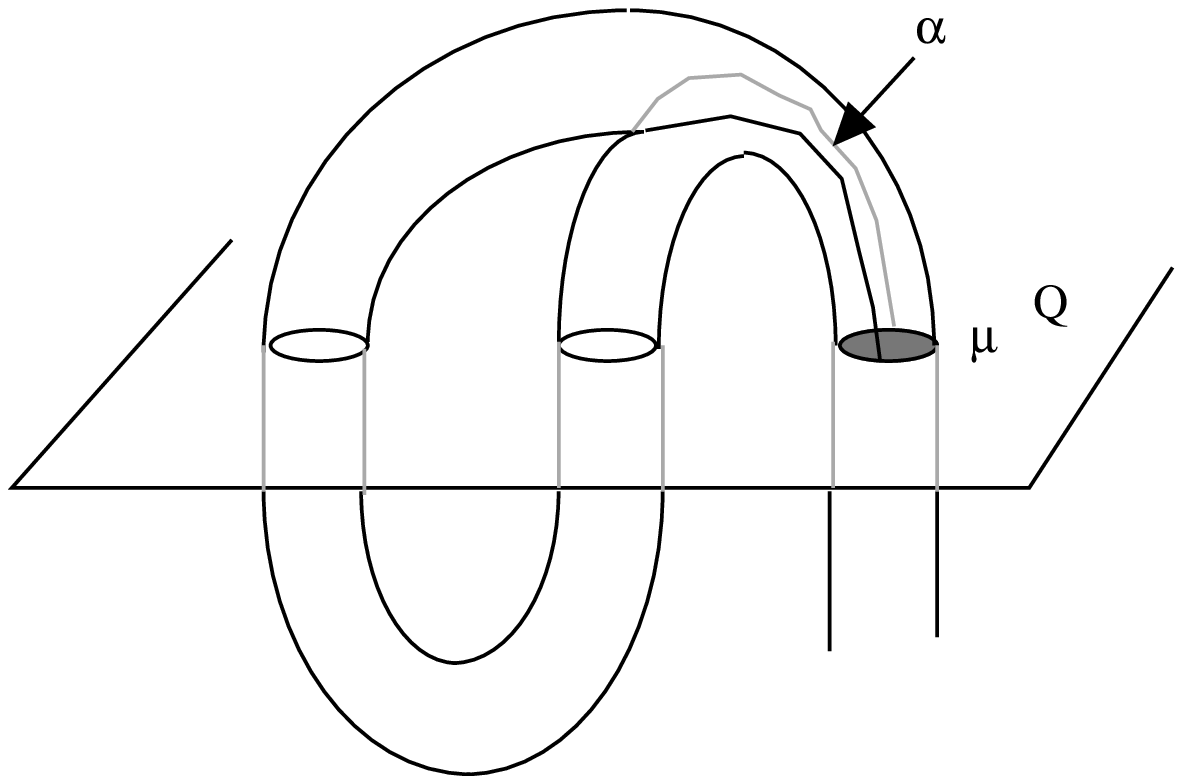}
\caption{wave}  \label{fig:wave}
\end{figure}

\begin{figure}
\centering
\includegraphics [width=0.65\textwidth]{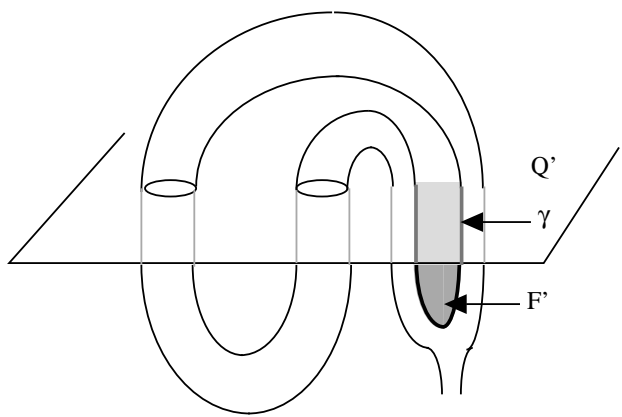}
\caption{Splitting along a wave}  \label{Q¹}
\end{figure}

Obviously $Q_{0} \cap e_{b}$ has one less point than $Q \cap e_{b}$. 
We claim that $Q_{0}$ is incompressible.  To verify this, consider the
tube dual to the compression disk (that is, $Q'$ is recovered by
tubing together two components of $Q_{0}$ along this tube).  The tube
is parallel to a regular arc $\ggg$ in $\bdd H$ connecting the two new
components of $Q_{0}$.  (The regular arc is one which intersects the
curve $\aaa$ in a single point.)  Let $F'$ be the disk of
parallelism, so $\bdd F'$ is the union of $\ggg$ and an arc in $Q'$
that crosses the compressing disk exactly once.  If there were a
further compression of $Q_{0}$ possible, it would have to fall on the
same side of $Q_{0}$ as $F'$.  Then note that $F'$ could be used to
push the compressing disk off the tube.  That is, the compression
could have been done to $Q$, which is impossible.  See Figure \ref{Q0}.

\begin{figure}
\centering
\includegraphics [width=0.65\textwidth]{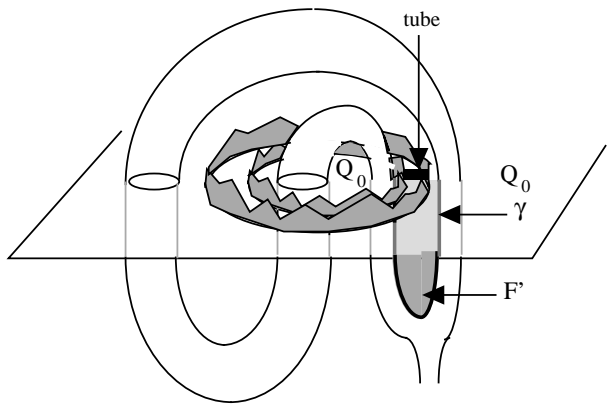}
\caption{Tubing the spheres $Q_{0}$ to get back $Q$}  \label{Q0}
\end{figure}

So the induction hypothesis applies to $Q_{0}$.  Since the first two
possibilities of the lemma imply (the first case of) the third, we may
as well take $F$ to be a disk as in the third possibility.  Note
specifically that if $Q_{0} \cap e_{b} = \emptyset$ then we can use
item 1) to choose for $F$ a disk that is disjoint from $e_{b}$.  When
comparing the curves $\aaa = F \cap H$ and $\ggg = F' \cap H$, we can
arrange that $\aaa \cap \ggg = \emptyset$ by pushing any intersection
points to the point where the tube is attached (to recover $Q'$ from
$Q_{0}$) and moving $\aaa$ across the attaching disks.  Note also that
at most one end of $\aaa \subset \bdd F$ lies on the new meridians of
$Q_{0}$ since these two meridians lie on different components of
$Q_{0}$.  By general position (make the tube thin) the interior of $F$
intersects $Q'$ only in meridians of the attaching tube.  Moreover,
since $F$ is disjoint from $\ggg$, all intersections of $F$ with $F'$
can be pushed via $F'$ across the tube so that, in the end, the
interior of $F$ is entirely disjoint from $Q'$ and from $F'$.  Now use
$F'$ to $\bdd$-compress $Q'$ to recover $Q$, leaving $F$ as a disk
satisfying the lemma for $Q$.
\end{proof}

We recall the definition of width for a graph; for further details see
\cite{GST}.  Let $\Ggg$ be an eyeglass graph or theta graph in
$S^{3}$.  As in \cite{GST}, choose a height function $h$ from $S^3$
with two points removed to $\mathbb{R}$, and let $S(t)=h^{-}(t)$.  Assume that
$\Ggg$ is in Morse position with respect to $h$, that is, the critical
points of $\Ggg$ with respect to $h$ occur at distinct values of $t$
and these values are distinct from the values of $h$ at the vertices
of $\Ggg$.  Further assume that a vertex $v$ of $\Ggg$ is either a
Y-vertex (where exactly two edges of $\Ggg$ lie above $v$) or a
$\lambda$-vertex (where exactly two edges of $\Ggg$ lie below $v$).

\begin{defin} 
Let $t_0<t_1<...<t_n$ be the successive critical heights of $\Ggg$ and
suppose $t_j$ and $t_k$ are the two levels at which the vertices
occur.  Let $s_i, 1\leq i\leq n$ be generic levels chosen so that
$t_{i-1}<s_i<t_i$.  Define the {\em width of $\Ggg$ with respect to h}
to be

$$
W_{h}(\Ggg)=2(\Sigma_{i\neq 
j,k}|S(s_i)\cap(\Ggg)|)+|S(s_j)\cap(\Ggg)|+ |S(s_k)\cap(\Ggg)|
$$

We say that $\Ggg$ is in {\em thin position} with respect to $h$ if 
has been isotoped to the generic position which minimizes $W_h$.
\end{defin}

\begin{figure}
\centering
\includegraphics [width=0.85\textwidth]{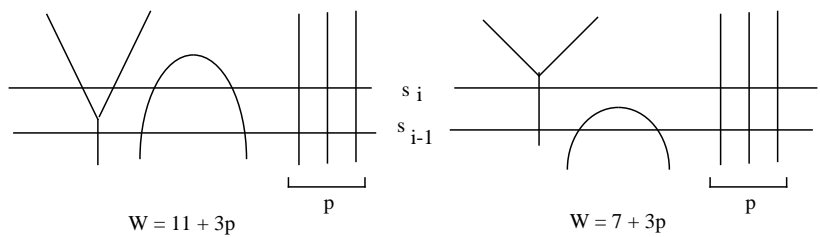}
\caption{Reducing the width by $4$ via Counting Rule
\ref{count:simple} case 4} \label{minmax}
\end{figure}

\begin{example}  If $\Ggg$ is a knot, then this definition of 
width is simply twice the width as defined by Gabai.
\end{example}

\begin{example} \label{exam:level} Suppose $e_{-}$ is a knot in
$S^{3}$, in generic position with respect to $h$.  Suppose $P$ is a
generic level sphere that intersects $e_{-}$ in $p$ points. 
Construct an eyeglass graph in $S^{3}$ by attaching to $e_{-}$ the
union of an edge $e_{b}$ and a loop $e_{+}$ both lying in $P$.  Then
when $\Ggg$ is made generic by tilting $e_{b} \cup e_{+}$,
$$W_{h}(\Ggg) = W_{h}(e_{-}) + 4p + 5.$$  Indeed, two vertices and a
regular maximum (say) are added.  Level spheres just below the
vertices add $p$ and $p+1$ to the width.  That just below the regular
maximum adds $2p + 4$.
\end{example}

We will mostly be concerned with how the width changes under 
isotopies of $\Ggg$, but it will be important to identify precise 
rules.  It is simple to check the following (see Figure \ref{minmax} for a 
sample argument):

\begin{num} \label{count:simple}  
    \begin{enumerate}
	
	\item As a maximum (either a regular maximum or a
	$\lll$-vertex) is pushed below (or above) another maximum, the
	width does not change.
	
	\item As a minimum (either a regular minimum or a
	$Y$-vertex) is pushed below (or above) another minimum, the
	width does not change.
	
	\item As a regular minimum is pushed above a regular maximum, the 
	width decreases by $8$.
	
	\item As a regular minimum is pushed above a $\lll$ vertex, or
	a regular maximum is pushed below a $Y$-vertex, the width
	decreases by $4$.
	
	\item As a $Y$ vertex is pushed above a $\lll$ vertex, or,
	equivalently, a $\lll$ vertex is pushed below a $Y$-vertex,
	the width decreases by $2$.  
	
	\item  Suppose between level spheres $P_{\pm}$ there are exactly two 
	critical points, a regular minimum and a regular maximum on the same 
	arc.  Then replacing that arc by a vertical arc reduces the width by 
	$4|P_{+} \cap \Ggg| + 4 = 4|P_{-} \cap \Ggg| + 4.$
	
    \end{enumerate}
    
\end{num}

\begin{defin} \label{defin:wequiv} Two embeddings of a trivalent graph
in $S^{3}$, both generic with respect to a height function on $S^{3}$,
are {\em width-equivalent} if there is a generic isotopy from one
embedding to the other so that the width is constant throughout the
isotopy.
\end{defin}

It's obvious that any birth-death singularity during the isotopy will
change the width, so the only non-generic embeddings during a
width-equivalence isotopy will be ones at which two critical points
are at the same level.  Note that, from Counting Principle
\ref{count:simple}, the two critical points must both be maxima or
both minima.  In other words, if two embeddings are width-equivalent
then there is an isotopy from one to the other that perhaps pushes
maxima past maxima and minima past minima, but never maxima past
minima.

\begin{defin} \label{defin:levellable} Suppose $\Ggg'$ is a subgraph
of a trivalent graph $\Ggg$ and $i_{1}: \Ggg \subset S^{3}$ is generic
with respect to the height function $h: S^{3} \map \mathbb{R}$.  We say that
$\Ggg'$ is {\em levellable} if there is an embedding $i_{2}: \Ggg \map
R^{3}$ so that
\begin{itemize}
    
    \item $i_{2}(\Ggg')$ is level.  That is, $hi_{2}(\Ggg') = t, t 
    \in  \mathbb{R}$
    
    \item $i_{1}$ is width-equivalent to an embedding obtained by 
    perturbing $i_{1}$
    
\end{itemize}
\end{defin}

For example, suppose $\Ggg \subset S^{3}$ is an eyeglass graph in
generic position with respect to $h$, except that one cycle $e_{\pm}$
in $\Ggg$ is level, e. g. $h(e_{+}) = t$.  There is a natural way to make
$\Ggg$ generic, namely tilt $e_{+}$ slightly so that it is vertical,
i. e. so that $e_{+}$ has a single maximum (perhaps a $\lll$ vertex) and
a single minimum (perhaps a $Y$-vertex) and one of these is the vertex
of $\Ggg$ lying in $e_{+}$.  The choice of whether the vertex is at
the minimum or at the maximum of $e_{+}$ is determined by whether the
end of the edge $e_{b}$ lies below or above the vertex.  The resulting
generic embedding of $\Ggg$ is one for which $e_{+}$ is levellable. 
In fact, using this convention, we can extend the notion of width so
that it is defined when either or both of $e_{\pm}$ are level.  An easy
application of Counting Rule \ref{count:simple} shows:

\begin{num} \label{count:horiz} Suppose that
$e_{+}$ is level and the end of $e_{b}$ at $e_{+}$ lies
below the vertex.
    \begin{enumerate}
	
	\item If $e_{+}$ is kept level while being moved below a
	regular maximum, the width increases by $4$.
	
	\item If $e_{+}$ is kept level while being moved below a
	$\lll$ vertex, the width increases by $2$.
	
	\item If $e_{+}$ is kept level while being moved above a
	regular minimum, the width increases by $8$.
	
	\item If $e_{+}$ is kept level while being moved above a
	$Y$ vertex, the width increases by $4$.
	
    \end{enumerate}
    
\end{num}

Of course the same rules apply when it is $e_{-}$ that is level, and
symmetric rules hold if the end of $e_{b}$ at the vertex lies above
the vertex.  There is one final case:

\begin{num} \label{count:wag} Suppose that $e_{+}$ is level, with
$h(e_{+}) = t$, and the end of $e_{b}$ at $e_{+}$ lies below the
vertex.  Let $p = |S(t) \cap (e_{-} \cup e_{b})|$.  If the end of
$e_{b}$ at $e_{+}$ is moved above $e_{+}$ (introducing a new regular
maximum in $e_{b}$) then the width is increased by $2p + 4$.
\end{num}

\begin{figure}
\centering
\includegraphics [width=0.85\textwidth]{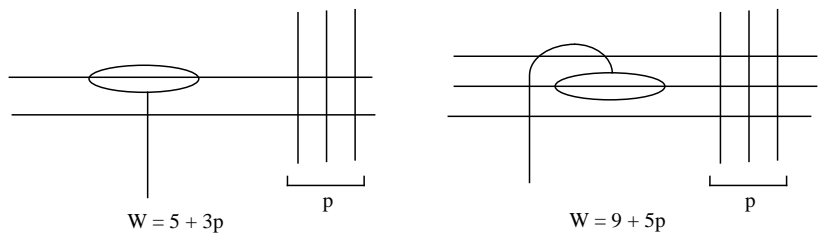}
\caption{Increasing $W$ by ``wagging'' the end of $e_{b}$}  \label{fig:wag}
\end{figure}

\begin{proof} See Figure \ref{fig:wag}
\end{proof}

\begin{lemma} \label{lemma:Haken}  
Let $\Ggg$ be a Heegaard spine eyeglass graph in $S^{3}$, in generic
position with respect to the height function $h$.  Suppose $e_{+}$
lies entirely above or entirely below $e_{-}$.  Then $\Ggg$ is planar
(i.e. can be isotoped to lie in a level sphere).
\end{lemma}

\begin{proof} The edge $e_{b}$
defines a Heegaard splitting of the reducible manifold $S^{3} -
\eta(e_{+} \cup e_{-}) \cong (S^{3} - \eta(e_{+}) \# (S^{3} - 
\eta(e_{+})$.  By Haken's theorem there is some reducing sphere that
intersects $e_{b}$ in a single point; planarity of $\Ggg$ (as well as 
the unknottedness of $e_{\pm}$) follows
immediately.
\end{proof}

\begin{lemma} \label{lemma:level} Let $\Ggg$ be the eyeglass graph
described in example \ref{exam:level} and suppose $\Ggg$ is a
Heegaard spine.  Suppose there is a maximum of $e_{-}$ below $P$ and
let $Q$ be a level sphere just above the highest such maximum. 
Suppose $|Q \cap \Ggg| = |Q \cap e_{-}| = q$.  Then the width of
$\Ggg$ can be reduced by at least $4q$.

(The symmetric statement hold if there is a minimum above $P$.)
\end{lemma}

\begin{proof} The proof is by induction on $q$.  Let $Q'$ be the 
collection of spheres obtained by maximally compressing $Q$ in the 
complement of $H$.  Note that $Q'$ intersects $\Ggg$ only in $e_{-}$, 
so Proposition \ref{prop:mori} case 1 applies.  The disk $F$ given by 
the proposition describes an isotopy that can be used to slide some 
part of an edge to $Q'$.  Hence, by avoiding the disks in $Q'$ that 
are the results of the compressions, the edge is brought down (or up) 
to $Q$.  The isotopy possibly passes through $Q$ on the way, but at 
the end $\bbb$ can be taken to lie just below (resp.  above) $Q$.  In 
particular, the arc moves down past a minimum (or at least past $e_{b} 
\cup e_{+}$) or it moves up past a maximum.  This decreases $q$ by $2$ 
and the width by at least $8$.  This would complete the inductive step 
unless the arc contains the maximum just below $Q$ (which would 
disrupt the induction) then the arc contains at least two minima as 
well as that maximum.  For the purposes of calculation of the 
resulting effect on width, we could imagine moving one of the 
contiguous minima up to just below the maximum (this will not thicken 
$\Ggg$) and then cancelling the minimum and maximum, thereby reducing 
the width by $4q + 4$, thereby accomplishing the required reduction.
\end{proof}

For a similar but more delicate argument that will soon follow we will
need to identify particularly important level spheres.

\begin{defin} Suppose $\Ggg \subset S^{3}$ is in generic position with
respect to the height function $h$.  A $Y-vertex$ at the minimum (or a
$\lll$-vertex at the maximum) of a vertical cycle is called an {\em
exceptional} critical point.  A generic level sphere $h^{-}(t)$ is
{\em thin} if the lowest critical point above it is a minimum and the
highest critical point below it is a maximum.  A thin level sphere is
{\em exceptional} if one (or both) of these critical points lying above or
below it is exceptional.
\end{defin}

\begin{lemma} \label{lemma:bridgeye}  
Let $\Ggg$ be a Heegaard spine eyeglass graph in $S^{3}$, in generic
position with respect to the height function $h$.  Suppose $e_{+}$ is
a vertical cycle with its minimum a $Y$-vertex $v$ and suppose that no
critical height of $\Ggg$ occurs between the heights of its minimum
and maximum.  Suppose there is some minimum of $\Ggg$ above $e_{+}$
and $P$ is the sphere just below the lowest such minimum.

Then either $\Ggg$ is planar or the width of $\Ggg$ can be reduced by
at least $2|\Ggg \cap P|$.

The symmetric statement is true for vertical cycles whose maximum is
a $\lll$-vertex.
\end{lemma}

\begin{proof} {\bf Special case:} $e_{-}$ is disjoint from $P$.

Following Lemma \ref{lemma:Haken} we may assume that $e_{-}$ does not 
lie entirely above $P$, so $e_{b}$ intersects $P$ in at least two 
points.  The proof in this case is by induction on $|\Ggg \cap P|$, 
and directly mimics the proof of Lemma \ref{lemma:level}.  Let $P'$ be 
the collection of spheres obtained by maximally compressing $P$ in the 
complement of $H$.  Since $P'$ intersects $\Ggg$ only in $e_{b}$, 
Proposition \ref{prop:mori} case 2 applies.  The disk $F$ given by the 
proposition describes an isotopy that can be used to slide some part 
of $e_{b}$ to $P'$.  Hence, by avoiding the disks in $P'$ that are the 
results of the compressions, the edge is brought down (or up) to $P$.  
In particular, the arc moves down past a minimum or it moves up past a 
maximum.  This decreases $|\Ggg \cap P|$ by $2$ and the width by at 
least $4$.  This completes the inductive step unless the arc contains 
the minimum just above $P$ (which would disrupt the induction).  But 
in this case the arc contains at least two maxima as well as that 
minimum.  For the purposes of calculation of the resulting effect on 
width, we could imagine moving one of the contiguous maxima down to 
just above the minimum (this will not thicken $\Ggg$) and then 
cancelling the minimum and maximum, thereby reducing the width by 
$4|\Ggg \cap P| + 4$, via Counting Rule \ref{count:simple} case 6, 
thereby accomplishing more than the required reduction, in this case.

So henceforth we assume that $e_{-}$ intersects $P$.  The structure of 
the argument will again mimic the proof of Lemma \ref{lemma:level}, 
though the details are a bit more complicated.  Let $Q_1, \ldots 
Q_{n}$, numbered from bottom to top, be the non-exceptional thin 
spheres for $\Ggg$.  That is, just above each $Q_i$ is a minimum that 
is not the $Y$-vertex minimum of a vertical cycle, and just below each 
$Q_{i}$ is a maximum that is not the $\lll$-vertex maximum of a 
vertical cycle.  So in particular $P$ is among these spheres.  Let $Q 
= Q_{1} \cup \ldots \cup Q_{n}$.  The proof will be by induction on 
$\Ggg \cap Q$.  Explicitly, we will show that given any 
counterexample, one can find a counterexample with fewer such 
intersection points.

Let $Q'$ be the collection of spheres obtained by maximally
compressing $Q$ in the complement of $H$.  Note that $Q'$ is disjoint
from $e_{+}$.  Let $F$ be the disk given by Proposition
\ref{prop:mori}.  There are two cases to consider:

{\bf Case 1:} $\aaa$ is a regular arc on $\bdd H$, disjoint from some 
meridian of $e_{b}$.

Then $F$ describes an isotopy that can be used to slide some part 
$e_{0}$ of an edge to $Q'$.  As usual, we can view this as bringing 
$e_{0}$ down (or up) to $Q$ so, at the end of the move, $e_{0}$ can be 
taken to lie just above (resp.  below) the $Q_{i}$ to which $e_{0}$ 
was adjacent.  In particular, the $e_{0}$ moves down past a minimum or 
up past a maximum.  If $\aaa = \bdd F \cap \Ggg$ does not go through a 
vertex (so $\aaa = e_{0}$), this reduces the width by at least $4$ 
($8$ if the critical point it passes is not a vertex) and it reduces 
$\Ggg \cap Q_{i}$ by $2$.  If $\aaa$ does pass through a vertex (so 
$e_{0} \subset \aaa$) the 
width drops by at least $2$ and $\Ggg \cap Q_{i}$ by $1$.  Note that 
$\aaa$ lies between $Q_{i}$ and one of $Q_{i \pm 1}$ so, unless $P = 
Q_{i}$ or $Q_{i-1}$, the move can have no effect on whether $P$ 
remains as described, or on $\Ggg \cap P$.  So unless $P = Q_{i}$ or 
$Q_{i-1}$ we are done, by induction.  In fact, even if $P = Q_{i}$ or 
$Q_{i-1}$ the result of the move gives a counterexample with $Q \cap 
\Ggg$ reduced, so long as $P$ remains as described.  That is, so long 
as a minimum remains just above $P$.

So suppose the slide or isotopy of $e_{0}$ to $\bbb \subset Q_{i}$ 
removes the last minimum above $P$ and suppose first that $P = 
Q_{i-1}$.  The effect is to remove $Q_{i-1}$ from $Q$ so the old 
$Q_{i}$ now serves as $P$.  We compute.  Let $p$ be the number of 
maxima between $Q_{i-1}$ and $Q_{i}$ before the move (counting any 
$\lll$ vertex as $1/2$ a maximum) and let $r$ be the number of minima 
(counting any $Y$ vertex as $1/2$ a minimum).  Then $|Q_{i-1} \cap 
\Ggg| - |Q_{i} \cap \Ggg| = 2p - 2r$.  We need to show that the move 
just described thins $\Ggg$ by at least twice that much, plus $4$ if 
$\aaa$ doesn't pass through a vertex (so $\Ggg \cap Q_{i}$ is reduced 
by two further points) or plus $2$ if $\aaa$ does pass through a 
vertex.  The computation is most obvious if $\aaa$ is a single minimum 
with both ends on $Q_{i}$, so $r = 1$ or $1/2$.  Then since this 
minimum passes $p$ maxima the width is reduced by at least $4p$ if the 
minimum is regular (even if the only maximum it passes is a 
$\lll$-vertex) and also $4p$ if the minimum is a $Y$-vertex, since we 
know that then all vertices are accounted for and the maxima are 
regular.  In any case, we have $4p \geq 2(2p - 2) + 4$, completing the 
computation in this case.

When $\aaa$ is more complicated, containing several minima, the only 
difference is an even greater thinning: for computational purposes one 
can imagine first moving a regular minimum in $e_{0}$ above all but 
its contiguous maxima, then cancelling the minimum with one of those 
contiguous maxima.  By Counting Rule \ref{count:simple} case 6, this 
already thins $\Ggg$ sufficiently; the actual isotopy would thin it 
even further.

The computation when $P = Q_{i}$ is similar.  In this case, if the
last minimum above $P$ is removed, $Q_{i+1}$ becomes the new $P$ and
we need to show that the width is reduced by at least $2(|Q_{i} \cap
\Ggg| - |Q_{i+1} \cap \Ggg|)$.  (We do not need to add $4$ or $2$,
since the move leaves $\Ggg \cap Q_{i+1}$ unchanged.)  If $Q_{i}$ was
the highest non-exceptional thin sphere then, for these computational
purposes, substitute a sphere above $\Ggg$ for $Q_{i+1}$.  Again let
$p$ and $r$ be the number of maxima and minima in the relevant region,
that is, between $Q_{i}$ and $Q_{i+1}$ (again, a $\lll$ vertex or $Y$
vertex counts as only half a maximum or minimum respectively.)  Since
the last minimum above $P$ is being eliminated by pulling $\aaa$ down
to $P$, a minimum of $\aaa$ has two contiguous maxima, which we may as
well take to be the highest two maxima between the spheres.  Then, for
computational purposes, we can imagine eliminating that minimum first,
dragging it past all but the two contiguous maxima, and then
cancelling it with one of the contiguous maxima.  The result is to
thin by at least $4p$ (in fact $8p$ if all relevant critical points
are regular) and this more than suffices.

{\bf Case 2:} $\aaa$ passes exactly once through a meridian of $e_{+}$ 
and has its ends at the same point of $e_{b} \cap Q_{i}$.

Then $P = Q_{i}$ or $Q_{i+1}$.  Suppose first that $P = Q_{i}$.  Then 
$F$ can be used to isotope the cycle $e_{+}$ so that it lies in 
$Q_{i}$, but now with the end of $e_{b}$ incident to it lying above 
$Q_{i}$.  When genericity is restored, $e_{+}$ is still vertical, but 
with its maximum now a $\lll$-vertex.  The simplest case to compute is 
when $\aaa$ runs through a single minimum of $e_{b}$, a minimum that 
lies just below $v$.  Then the move described eliminates that regular 
minimum, so one less term appears in the calculation of width.  This 
is the reverse of the operation described in Counting Rule 
\ref{count:wag}, so the width is decreased by $2|Q_{i} \cap \Ggg| + 4 
= 2|P \cap \Ggg| + 4$, immediately confirming the lemma.  If the end 
of $e_{b}$ near $e_{+}$ is more complicated, the thinning is even 
greater.

\begin{figure}
\centering
\includegraphics [width=0.65\textwidth]{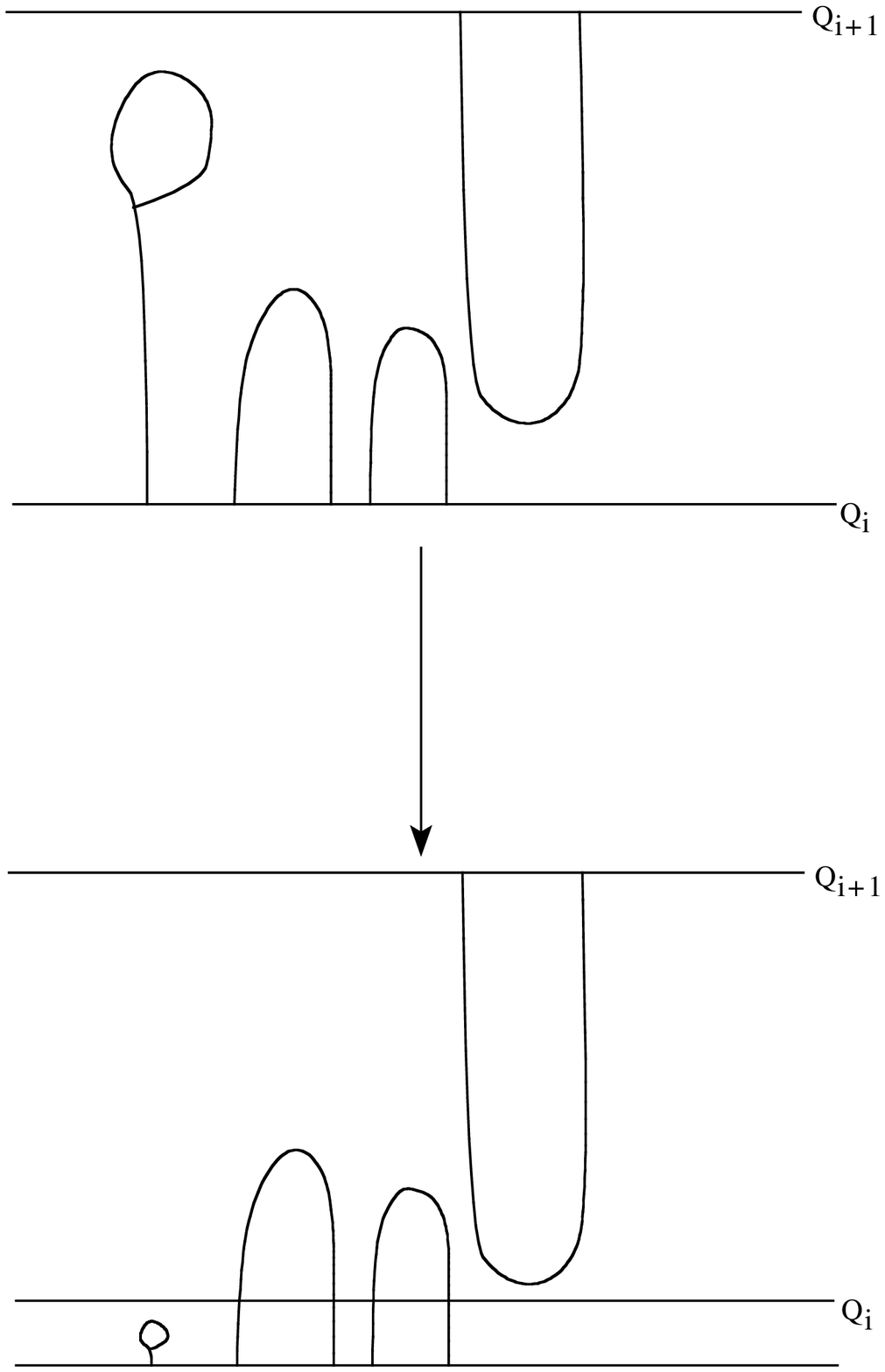}
\caption{}  \label{fig:thm6}
\end{figure}

Finally suppose that $P = Q_{i+1}$.  In this case the isotopy given by 
$F$ pulls $e_{+}$ down to $Q_{i}$.  Again consider the simplest case: 
the end segment of $e_{b}$ eliminated by the move is a simple vertical 
arc between $Q_{i}$ and $e_{+}$.  Then $F$ pulls $e_{+}$ past $p$ 
maxima and $r$ minima, changing the width by $4p - 8r$, essentially by 
Counting Rule \ref{count:horiz}.  (Again a $\lll$ vertex and $Y$ 
vertex count as only half a maximum or minimum respectively.)  On the 
other hand, $|Q_{i} \cap \Ggg|$ differs from $|Q_{i+1} \cap \Ggg|$ by 
$2p - 2r$ and $4p - 8r < 2(2p - 2r)$.  So, after the move, we have an 
even more extreme counterexample, and one with fewer points of 
intersection with $Q$.  Furthermore, if $e_{b}$ is in fact more 
complicated than a simple vertical arc, then even more thinning would 
have been done.  Now apply the inductive hypothesis and the 
contradiction completes the proof.
\end{proof}

\section{Main Theorems}

\begin{theorem} \label{theorem:bridgeye}  
Let $\Ggg$ be a tri-valent graph that is a genus two Heegaard spine in
$S^{3}$.  If $\Ggg$ is in thin position then it is in extended bridge
position.
\end{theorem}

\begin{proof} Suppose $\Ggg$ is not in extended bridge position.  As
previously, let $Q_1, \ldots Q_{n}$, numbered from bottom to top, be
the non-exceptional thin spheres and let $H = \eta(\Ggg)$.

Suppose some $Q_{i}$ intersects $H$ only in non-separating meridians.  
Then the argument is much as in the Special Case of Lemma 
\ref{lemma:bridgeye}: Let $Q'$ be the collection of spheres obtained 
by maximally compressing $Q_{i}$ in the complement of $H$.  By 
Proposition \ref{prop:mori} each component of $Q'$ is parallel to a 
component of $\bdd H - Q'$.  So in particular, there is a disk $F 
\subset S^{3}$ whose interior is disjoint from $H \cup Q'$ and $\bdd F 
= \aaa \cup \bbb$, where $\aaa \subset \bdd H$, $\bbb \subset Q'$ and 
$\aaa$ is a regular path on $\bdd H$ (not intersecting some meridian 
of $e_{b}$, if $\Ggg$ is an eyeglass).  Since $\aaa$ is disjoint from 
$\bdd Q' = \bdd Q_{i}$, $\aaa$ lies entirely above or below, say 
above, the level of $Q_{i}$.  Then $F$ describes an isotopy that can 
be used to slide some part $e_{0}$ of an edge down to $Q_{i}$.  The isotopy 
possibly passes through $Q_{i}$ on the way, but at the end $e_{0}$ can 
be taken to lie just above $Q_{i}$.  In particular, $e_{0}$ either 
lies below the minimum just above $Q_{i}$ or the arc containing that 
minimum has been changed to one with a single maximum just above 
$Q_{i}$.  In any case, the graph is thinned, a contradiction.

So assume every $Q_{i}$ intersects $H$ in some separating meridians,
that is, $\Ggg$ is an eyeglass graph and for each $i$, $Q_{i} \cap
e_{b} \neq \emptyset$.

If any $Q_{i}$ is disjoint from both of $e_{\pm}$, we use the
same argument as in the Special Case of Lemma \ref{lemma:bridgeye}, 
with $Q_{i}$ playing the role of $Q$.

So assume every $Q_{i}$ intersects either $e_{+}$ or $e_{-}$ as well 
as $e_{b}$.  If each $e_{\pm}$ intersects some $Q_{i}$, we use the 
same argument as above, via Proposition \ref{prop:mori} case 3.  We 
are left with the case that $e_{+}$, say, is disjoint from all 
$Q_{i}$, whereas $e_{-}$ intersects every $Q_{i}$.  So suppose $e_{+}$ 
lies between $Q_{i}$ and $Q_{i+1}$ and, for concreteness and with no 
loss of generality (by symmetry) assume that the point $q$ of $Q \cap 
e_{b}$ that is closest to $e_{+}$ lies in $Q_{i}$, some $1 \leq i \leq n$.  
(Here if $i = n$, $Q_{i + 1}$ is taken to be a level sphere above 
$\Ggg$.)

\bigskip

{\bf Claim:} $e_{+}$ is a vertical cycle lying above some maximum of
$\Ggg$ that lies between $Q_{i}$ and $Q_{i+1}$.  The minimum of
$e_{+}$ is a $Y$-vertex.

\bigskip

{\bf Proof of claim} Let $Q'$, as before, be the collection of spheres 
obtained by maximally compressing $Q_{i}$ in the complement of $H$.  
As we have argued, Proposition \ref{prop:mori} shows that there is a 
disk $F$ for $Q'$ as given in item 3b of that Proposition.  That is, 
$\bdd F$ consists of an arc $\bbb$ on $Q_{i}$ with both ends at $q$ 
and an arc $\aaa$ on $\bdd H - Q$ parallel to a cycle with both ends 
at $q$ and running once around $e_{+}$.  $F$ can be used to pull the 
component of $\Gamma - Q_{i}$ that contains $e_{+}$ down to $Q_{i}$.  
For computational purposes we can picture this done in three stages: 
$e_{+}$ is replaced by a vertical cycle with its minimum (resp.  
maximum) at the minimum (resp.  maximum) of $e_{+}$; the end of 
$e_{b}$ between $Q_{i}$ and $e_{+}$ is replaced by a vertical arc 
terminating at the minimum of $e_{+}$; and then $e_{+}$ and the end of 
$e_{b}$ are pulled down to $Q_{i}$.  The first two steps cannot make 
$\Ggg$ thicker and will make it thinner unless in fact it leaves the 
height function on $\Ggg$ unchanged.  The third move will not thicken 
$\Ggg$ if the original $e_{+}$ has a minimum below all the maxima (e.  
g.  there is a regular minimum of $e_{+}$) and in fact must thin 
$\Ggg$ unless $e_{+}$ lies above some maximum.  So, since $\Ggg$ 
cannot be thinned, $e_{+}$ must be a cycle containing no regular 
minima and lying entirely above some maximum.  This proves the claim.

\smallskip

Having established the claim, Lemma \ref{lemma:bridgeye} applied to $P 
= Q_{i+1}$ implies that $i = n$ so $\Ggg \cap P = \emptyset$.  But 
even then, the {\em argument} of Lemma \ref{lemma:bridgeye} still 
suffices to display the same contradiction: The effect of pulling 
$e_{+}$ to $Q_{i}$ is to alter the width by adding at most $4p - 8r$.  
On the other hand, after the move, $Q_{i}$ is then suitable (when 
pushed just above $e_{+}$) for applying Lemma \ref{lemma:bridgeye}.  
(See Figure \ref{fig:thm6}.)  This lemma says that $\Ggg$ can be 
thinned by $2|Q_{i} \cap \Ggg| = 4p - 4r > 4p - 8r$.
\end{proof}

\begin{defin} Suppose $\Ggg$ is in bridge position.  Then a
level sphere separating the minima from the
maxima is called a {\em dividing sphere} for $\Ggg$.  

If $\Ggg$ is not in bridge position, but is in 
extended bridge position, then a dividing sphere is a level sphere $P$ for 
which every minimum above $P$ is the $Y$-vertex of a vertical cycle 
and every maximum below $P$ is the $\lll$-vertex of a vertical cycle. 
\end{defin}

\begin{theorem} \label{theorem:thineye}  
Let $\Ggg$ be a tri-valent graph that is a genus two Heegaard spine in
$S^{3}$.  If $\Ggg$ is in thin position then it is in extended bridge
position.  Either $\Ggg$ is planar or some dividing sphere is disjoint
from a simple (i.  e. non-loop) edge of $\Ggg$.
\end{theorem}

\begin{proof} Following Theorem \ref{theorem:bridgeye} we can assume
that $\Ggg$ is in extended bridge position.  If $\Ggg$ is in
(non-extended) bridge position, the proof (and Corollary
\ref{cor:leveledge}) will conclude much as in Theorems \cite[5.3,
5.14]{GST}.  We note that were we content to find {\em either} a 
level edge {\em or} an unknotted cycle in $\Ggg$, we would be done
following this case.  However the pursuit of a simple edge requires
more persistence.  Since the delicate points in the argument will need
to be repeated in the case of extended bridge position we only
summarize the proof when $\Ggg$ is in bridge position:

There is a dividing sphere $Q$ between the lowest maximum and the
highest minimum that cuts off both an upper disk and a lower disk.  If
an edge running between distinct vertices lies above or below $Q$ we
are done.  So we can assume that each component of $\Ggg - Q$ is
either an arc or a $3$-prong.  (This fact makes some of the
complications in the proof of \cite[5.3]{GST} irrelevant.)  There is
an argument to show that we can find such upper and lower disks so
that their interiors are disjoint from $Q$ and that neither intersects
$Q$ in a loop.  Each is incident to exactly two points of $\Ggg \cap
Q$ and it is shown that at least one point, and perhaps both, are the
same for both upper and lower disks.

If both upper and lower disks are incident to the same pair of points,
then these disks can be used to make a cycle (either a loop or a
$2$-cycle) level.  The argument of \cite[5.14]{GST} shows that if the
cycle is a loop then either $\Ggg$ could be thinned (a contradiction
to hypothesis) or $e_{b}$ is already disjoint from the dividing sphere
and we are done.  Essentially the same argument applies in the case of
a level $2$-cycle, unless the third edge too can be moved into the
sphere.  In the latter case, the graph is planar.

If the upper and lower disks are incident to only one point of $\Ggg
\cap Q$ in common, then they may be used either to thin $\Ggg$ or to
make that edge level, lying in $Q$.  In this case, too, $\Ggg$ may be
thinned, or another edge brought to $Q$ (creating a level $2$-cycle)
this time by using an outermost disk of a meridian $E$ for $S^3 - H$,
cut off of $E$ by $Q - H$.  For details see \cite[6.1, Subcases 3a, 3b]{GST}.

\begin{figure}
\centering
\includegraphics [width=0.65\textwidth]{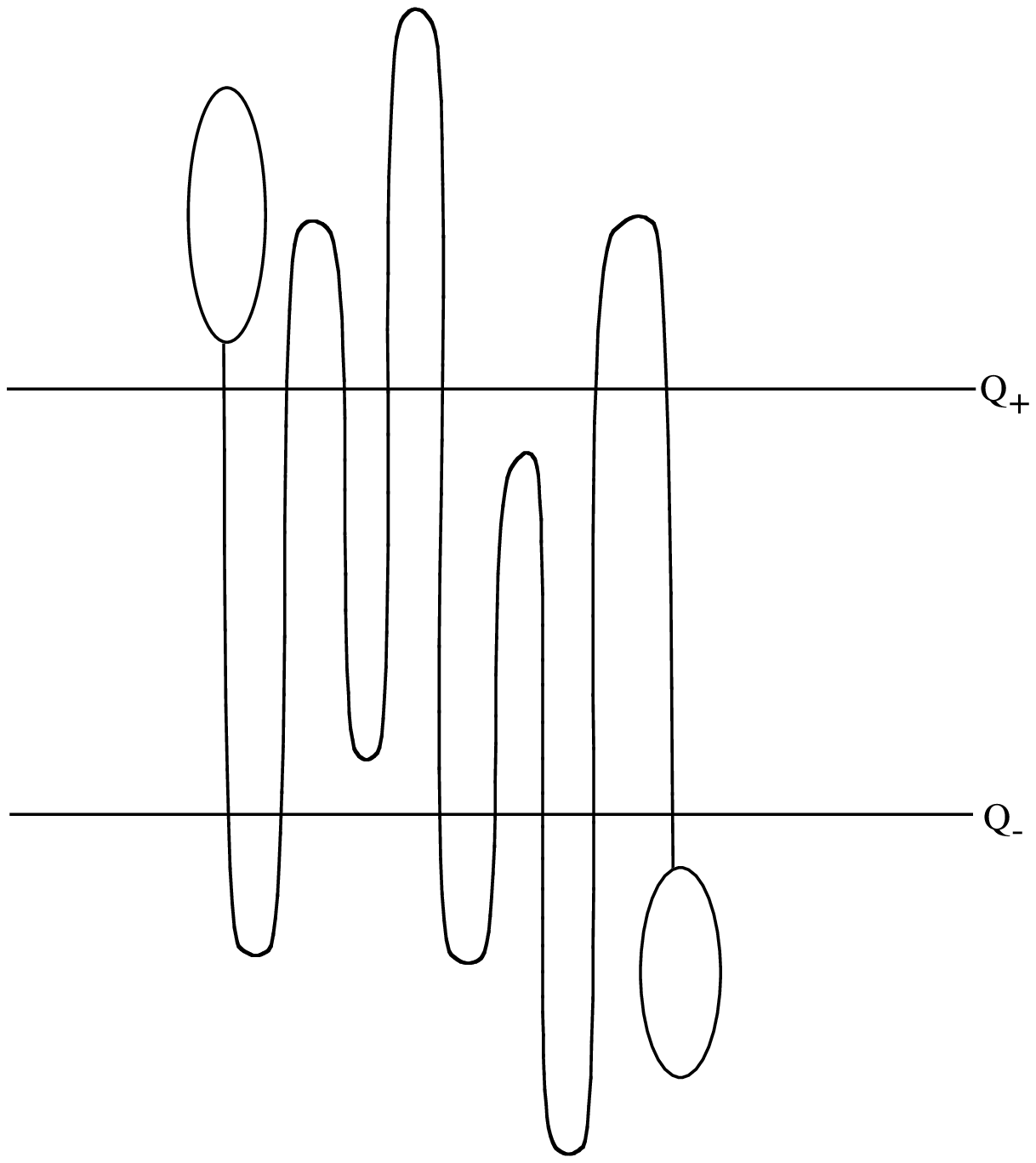}
\caption{}  \label{thm7.1}
\end{figure}

So now assume that $\Ggg$ is not in bridge position, but only in
extended bridge position.  In particular, all thin spheres are
exceptional and there is at least one exceptional thin sphere.

\bigskip

{\bf Claim 1:}  There is exactly one exceptional thin sphere and it 
intersects exactly one of the loops $e_{\pm}$.  

{\bf Proof of Claim 1:} Since there are at most two vertical cycles,
there are at most two exceptional thin spheres.  If there are two,
denote them by $Q_{\pm}$, with $Q_{+}$ lying above $Q_{-}$ (Figure
\ref{thm7.1}).  Consider the lowest minimum above $Q_{+}$ and the
highest maximum below $Q_{-}$.  It can't be that neither of these
critical point is exceptional, for then $\Ggg$ would not be in
extended bridge position.  If both critical points are exceptional,
then $\Ggg$ is planar by Lemma \ref{lemma:Haken}.  So we may as well
assume that both exceptional vertices are exceptional minima, one just
above $Q_{-}$ and one just above $Q_{+}$.  But then $Q_{-}$ intersects
$\Ggg$ only in $e_{b}$, contradicting thin position, via Proposition
\ref{prop:mori} case 2.  

Having established that there is exactly one 
exceptional thin sphere, the same argument shows that it cannot be 
disjoint from both $e_{\pm}$.  

\bigskip

With no loss of generality, suppose $e_{+}$ but not $e_{-}$ is 
disjoint from the exceptional thin sphere $Q$.

\bigskip

{\bf Claim 2:} The loop $e_{+}$ can be isotoped to lie in 
$Q$, without increasing the width of $\Ggg$.

{\bf Proof of Claim 2:} Maximally compress the exceptional level
sphere $Q$ in the complement of $H$ and call the result $Q'$.  Apply
Lemma \ref{prop:mori} to deduce that there is a disk $F$ as in item 3. 
Since it cannot describe a way to slide an edge segment of $\Ggg - Q$
to the level of $Q$ (that would make $\Ggg$ thinner), $\bdd(F)$ must
be disjoint from $e_{-}$ and run around $e_{+}$.  $F$ can then be used
to isotop $e_{+}$, as required.  Since the vertex of the loop is
immediately adjacent to $Q$, this does not thicken $\Ggg$.

\bigskip

Following the isotopy of Claim 2, $e_{+}$ divides $Q$ into two disks,
$Q_{1}$ and $Q_{2}$.  Consider the intersection of these $Q_{i}$ with
a meridian disk $E$ of $S^{3} - H$.  Note that there can be no closed
components of intersection, since an innermost one, if essential in
$Q_{i} - \Ggg$, could be used to push part of $\Ggg$ through $Q_{i}$,
thinning $\Ggg$.  (It is thinned, per Counting Rule \ref{count:horiz},
because an upper cap would contain no minima, and a lower cap would
contain more minima than maxima).  Similarly, an outermost arc of $E -
Q_{i}$ can't cut off a disk lying entirely above $Q$, for it could be
used to thin and, indeed, unless $Q$ and $e_{b}$ are disjoint, so
could a lower one, essentially by Counting Rule \ref{count:wag} 
applied in reverse.
 
\begin{figure}
\centering
\includegraphics [width=0.65\textwidth]{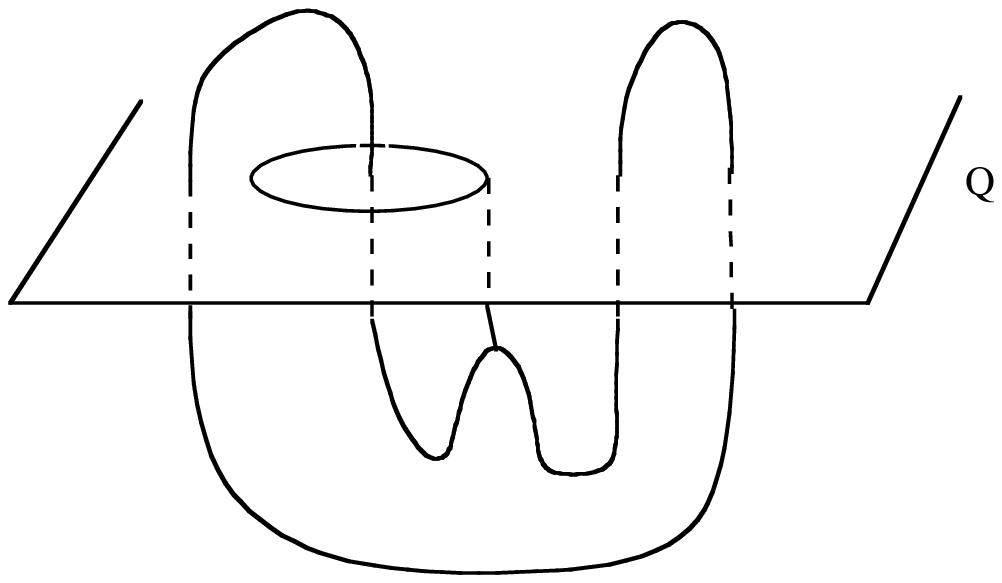} 
\caption{}  \label{thm7.3}
\end{figure}

So we may as well assume that $e_{b} \cap Q = \emptyset$.  We know 
that a maximum lies just below $Q$.  One possibility is that there is 
a regular maximum below $Q$.  Another is that the only maximum below 
$Q$ is a $\lll$-vertex (Figure \ref{thm7.3}).  

In the second case, if the end of $e_{b}$ is incident to the
$\lll$-vertex from above, then $e_{b}$ is monotonic (for otherwise an
internal maximum would lie below $Q$ or $e_{b}$ would intersect $Q$,
both possibilities we are not considering).  Then $e_{b}$ is disjoint
from the level sphere (a dividing sphere) just below the
$\lll$-vertex, and we are done.  So either there is a regular maximum
below $Q$ or the $\lll$-vertex below $Q$ has the end of $e_{b}$
incident to the vertex from below.  In particular, a level sphere just
below either sort of maximum would cut off an upper disk.  So, as is
now standard, some dividing sphere $P$ can be placed so that it
simultaneously cuts off both an upper disk $D_{u}$ and a lower disk
$D_{l}$.  As noted above, we can assume that neither disk has a closed
curve of intersection with $Q$.  We now proceed to duplicate, in this
context, the proof of \cite[5.3]{GST}.  The added difficulties here
are apparent even at the first step.  We will consider the
intersections of the interiors of $D_{u}$ and $D_{l}$ with $P$.

\bigskip

{\bf Claim 3:} (cf.  \cite[Claim 5.5]{GST}) There cannot be both an 
upper cap and a lower cap whose boundaries are disjoint.

{\bf Proof of Claim 3:} Let $C_{u}$ and $C_{l}$ denote the caps.  They
bound disjoint disks $P_{u}$ and $P_{l}$ in $P$.  If the end segment
of $e_{b}$ at $e_{+}$ is not incident to $P_{u}$ the proof is natural:
pushing $C_{u}$ down to $P_{u}$ and $C_{l}$ up to $P_{l}$ will thin
$\Ggg$.  So assume that $e_{+}$ does lie between $C_{u}$ and $P_{u}$. 
If any maximum is incident to $P_{u}$ and is lower than the height of
$Q$ (i.  e. the height of $e_{+}$) then $\Ggg$ could be thinned by
just pushing that maximum down while pushing $C_{l}$ up.  So any
maximum lying between $C_{u}$ and $P_{u}$ is higher than $e_{+}$.  On
the other hand, if any maximum {\em not} between $C_{u}$ and $P_{u}$
were above $e_{+}$ it could be pushed lower (since its easy to make
the descending disk from that maximum disjoint from $C_{u}$.  This too
would thin $\Ggg$.  Hence we see that the $p \geq 0$ maxima that are
lower than $e_{+}$ are precisely those that don't lie between $C_{u}$
and $P_{u}$.

Now consider the effect of pushing $C_{u}$ down to $P_{u}$ while
simultaneously pushing $C_{l}$ up to $P_{l}$.  Apply Counting Rule
\ref{count:horiz}: Pushing $e_{+}$ past $p$ maxima increases the width
by $4p$ whereas pushing up the $r \geq 1$ minima between $C_{l}$ and
$P_{l}$ reduces the width by $8r$.  (Here, as was usual in such
counting above, a $\lll$-vertex or $Y$ vertex counts as only half a
maximum or minimum).  The result is that, after the push, the width is
increased by at most $4p - 8r$.  On the other hand, after the push,
$P$ would satisfy the hypotheses of Lemma \ref{lemma:bridgeye}.  It's
easy to calculate $P \cap \Ggg$: it's $2p -2r$.  Then according to
that lemma, $\Ggg$ could be thinned by a further $4p - 4r > 4p -8r$, a
contradiction establishing the claim.

\bigskip

{\bf Claim 4:} (cf.  \cite[Claim 5.6]{GST}) If there is an upper disk 
and a disjoint lower cap, then we can find such a pair for which the 
interior of the upper disk is disjoint from $P$.  (The symmetric 
statement is of course also true.)

{\bf Proof of Claim 4:} Let $B_{u}$ and $B_{l}$ denote the balls above
and below the dividing sphere $P$ respectively.  The proof would
follow just as in \cite{GST} if we could find a complete collection
$\Ddd$ of descending disks for $\Ggg \cap B_{u}$ such that the
boundaries of $\Ddd$ and $D_{u}$ intersect only on $P$.  We do not
need to worry here, as we did there, about components of $\Ggg \cap
B_{u}$ that contain two vertices for if such a component exists the
lemma is proven.  What we do need to worry about is that any maxima
that are higher than the loop $e_{+}$ have no descending disks at all
(or rather, their descending disks encounter $e_{+}$ at $Q$ and do not
descend to $P$, else we could thin $\Ggg$.)  But because we have
established above that $D_{u}$ is disjoint from $Q$ there is an easy
fix.  The graph $\Ggg$ intersects the region $S_{P}^{Q} \cong S^{2}
\times I$ between $Q$ and $P$ in a collection of maxima and a
collection of vertical arcs.  At the top of one vertical arc
$\epsilon_{b}$ (an end of $e_{b}$) we see the bottom half of the loop
$e_{+}$.  Let $T$ be the union of two trees in $Q - e_{+}$, each
having a root at the vertex in $e_{+}$, each on opposite sides of
$e_{+}$ and together containing all the other points of $\Ggg \cap Q$. 
(These points are just the tops of the vertical arcs of $\Ggg \cap
S_{P}^{Q}$.)  Denote the edges of $T$ by $E_{T}$.  Finally, let $C
\subset S_{P}^{Q}$ be the vertical cylinder $e_{+} \times I$,
intersecting $\Ggg$ exactly in $\epsilon_{b} \cup e_{+}$.  Define
$\Ddd$ to be this collection of disks: $\{ E_{T} \} \times I$, $C -
\eta(\epsilon_{b})$, and a set of descending disks for all maxima in
$S_{P}^{Q}$, these chosen to be disjoint from the other disks in
$\Ddd$.  Clearly $\Ddd$ cuts $S_{P}^{Q}$ up into a collection of
balls.  See Figure \ref{fig:claim4}.

\begin{figure}
\centering
\includegraphics [width=0.65\textwidth]{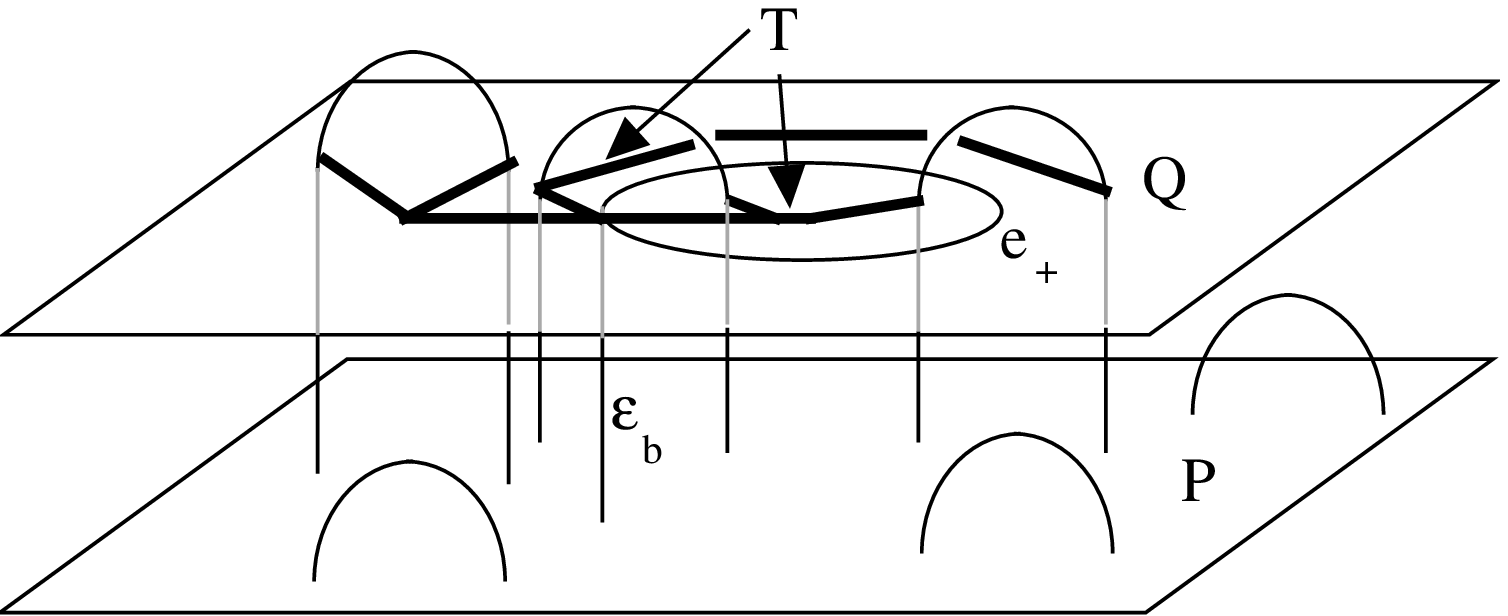} 
\caption{}  \label{fig:claim4}
\end{figure}

Now observe that $D_{u}$ cannot involve the maxima that are higher
than $Q$, else $\Ggg$ could be thinned.  Hence the part of the
boundary of $D_{u}$ that lies on $\Ggg$ either lies on a maximum in
$S_{P}^{Q}$ or on the component containing $e_{+}$.  In either case it
is easily made disjoint from $\bdd \Ddd$ so that $\bdd D_{u} \cap \bdd
\Ddd$ lies entirely in $P$.  The proof now follows as in \cite[Claim
5.6]{GST}.

\bigskip 

With one exception, the proof of Theorem \ref{theorem:thineye} is now
little different from the flow of the proof of \cite[Theorem
5.3]{GST}: ultimately we get upper and lower disks which can be used
to push part of $\Ggg \cap B_{l}$ up while pushing part of $\Ggg \cap
B_{u}$ down.  Unless the latter is the component containing $e_{+}$,
this immediately thins $\Ggg$.  So suppose $D_{u}$ does push down
$e_{+}$; let $p$ denote the point in $e_{b}$ where that component is
cut off.  Unless $D_{l}$ pushes up a segment incident to $p$, the
proof follows by a width count and Lemma \ref{lemma:bridgeye} just as
in the proof of Claim 3.  If the segment incident to $p$ that $D_{l}$
pushes up is a simple minimum (i.  e. it does not contain the other
end of $e_{b}$) then that push eliminates a critical point which we
may take to lie just below $P$.  In particular, for $P_{+}$ a level
sphere just above $P$, the move reduces the width by $2|\Ggg \cap
P_{+}| + 4$ via Counting Rule \ref{count:horiz}, and this is enough
again to ensure that after the move the graph is thinner. 

Finally, suppose that $D_{l}$ is incident to $p$ and pushes up the
other end of $e_{b}$.  (This implies in particular that $e_{b} \cap P
= \{ p \}$.)  Then after the move both the edges $e_{b}\cup e_{+}$ are
level and lie in $P$.  But, as usual, the move may thicken $\Ggg$ and
this time there is no immediate cancellation of a critical point since
$e_{b}$ was monotonic before the move, just as it would be again when
genericity is restored.  The thickening occurs, as usual, because the
$Y$-vertex minimum of $e_{+}$ may be pulled down past $m$ maxima, in
which case the width increases by $4m$.  But, unless $m = 0$, this
leads to a contradiction: Consider the cylinder $C$ that is swept out
by $e_{+}$ as it is pulled down to $P$ (effectively, this is just
another way of viewing the upper disk $D_{u}$) and apply the technical
Lemma \ref{lemma:tech} that follows.  The resulting graph could in fact
be thinned by a further $4m + 4$, leaving it thinner than we started. 
So we conclude that $m = 0$ and the move can be made without any
thickening at all.

\begin{figure}
\centering
\includegraphics [width=0.85\textwidth]{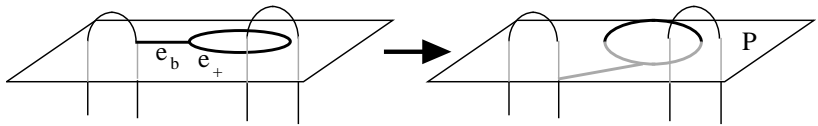} 
\caption{Tilting $e_b \cup e_{+}$}  \label{fig:newP}
\end{figure}

Once $e_{b} \cup e_{+}$ is level, tilt it slightly, creating two 
$Y$-vertices, say, one at each end of $e_{b}$, so $e_{+}$ is vertical 
with its maximum a regular maximum.  Then a level sphere passing 
through the middle of $e_{+}$ is a dividing sphere that is disjoint 
from $e_{b}$, as required.  See Figure \ref{fig:newP}.
\end{proof}

\begin{figure}
\centering
\includegraphics [width=0.85\textwidth]{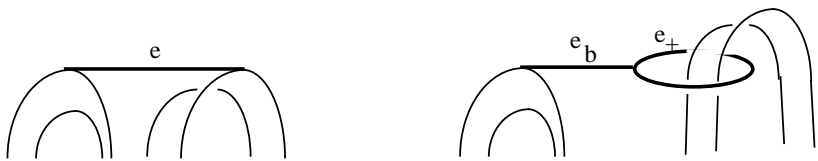} 
\caption{Levellable edge or subgraph}  \label{fig:leveledge}
\end{figure}

\begin{cor}\label{cor:leveledge}
Let $\Ggg$ be a tri-valent graph that is a Heegaard spine in $S^{3}$
and suppose that $\Ggg$ is in thin position.  Then at least one
simple edge is levellable (cf Definition \ref{defin:levellable}). 
To be specific, either $\Ggg$ is planar or (see Figure
\ref{fig:leveledge}):

\begin{enumerate}

\item If $\Ggg$ is in bridge position then there is a simple edge $e
\subset \Ggg$ so that

\begin{itemize}
    \item the knot or link $K = \Ggg - interior(e)$ is in bridge position and
    
    \item $e$ is levellable and its ends lie at distinct maxima or at
    distinct minima of $K$
    
\end{itemize}
    
\item If $\Ggg$ is not in bridge position then $\Ggg$ is an eyeglass graph.  
For some loop (say $e_{-}$) in $\Ggg$ 

\begin{itemize}
    \item $e_{-}$ is in bridge position and
    
    \item the subgraph $e_{b} \cup e_{+}$ is levellable and is
    incident to either a maximum or minimum of $e_{-}$.
    
\end{itemize}
\end{enumerate}
\end{cor}

\begin{proof} We assume $\Ggg$ is not planar and first suppose $\Ggg$
is in bridge position.  Let $P$ be a dividing sphere disjoint from a
non-loop edge $e$ of $\Ggg$ guaranteed by Theorem
\ref{theorem:thineye}.  With no loss of generality the edge $e$ lies
above $P$.  Let $\Ggg_{u}$ denote the part of $\Ggg$ lying above $P$. 
Since there are no minima above $P$, a family of descending disks for
$\Ggg_{u}$ describes a parallelism between $\Ggg_{u}$ and a subgraph
of $P$.  In particular, $e$ can be viewed as a perturbed level edge.

Suppose next that $\Ggg$ is not in bridge position.  We know from
Theorem \ref{theorem:bridgeye} that $\Ggg$ is extended bridge position
so in particular $\Ggg$ is an eyeglass graph.  Let $P$ be a dividing
sphere disjoint from the edge $e_{b}$, as guaranteed by Theorem
\ref{theorem:thineye}.  We may as well assume $e_{b}$ lies above $P$,
so one end of $e_{b}$ descends from the minimum of a vertical loop,
say $e_{+}$.  Since $e_{b}$ is disjoint from the dividing sphere $P$
it contains no minimum and its other end ascends from a $\lll$-vertex,
hence from a maximum of $e_{-}$.  Raise that maximum along $e_{b}$
until it is the critical point just below the $Y$-vertex.  Let $Q$ be
a level plane that intersects the monotonic edge $e_{b}$ in a single
point.  Maximally compress $Q$ in the complement of $\Ggg$ and let the
result be $Q'$.  As has been argued repeatedly above, if we apply
Proposition \ref{prop:mori} to $Q'$ the only conclusion that does not
violate thinness is possibility 3.b. In that case, the disk $F$
describes how to isotope $e_{b} \cup e_{+}$ to lie in $P$.  Since
there are no critical points between the heights of the ends of
$e_{b}$ this has no effect on width.
\end{proof}

\section{Technical Lemma}

For the following technical lemma we return to the context of Example
\ref{exam:level} and Lemma \ref{lemma:level}.  That is, $e_{-}$ is
generic with respect to a height function on $S^{3}$ and the subgraph
$e_{b} \cup e_{+}$ is level with respect to the height function, at a
height that is generic for $e_{-}$.  Width is calculated by tilting
$e_{b} \cup e_{+}$ slightly to restore genericity.  This is
independent of the direction of tilting.

\begin{lemma} \label{lemma:tech}
Suppose $\Ggg$ is a non-planar eyeglass graph that is a Heegaard spine
of $S^{3}$.  Suppose there is a height function on $S^{3}$ and a
dividing sphere $P$ for $e_{-}$ that contains both the edges $e_{b}$
and $e_{+}$.  Suppose $Q$ is a level sphere above $P$ and there is a
properly embedded annulus $C$ such that

\begin{enumerate}

\item $C$ spans the region $S_{P}^{Q} \cong S^{2} \times 
I$ that lies between $Q$ and $P$

\item $\bdd C \cap P = e_{+}$ and

\item $C \cap \Ggg = e_{+}$.  

\end{enumerate}

Let $m > 0$ be the number of maxima of $e_{-}$ in $S_{P}^{Q}$.  Then 
$\Ggg$ can be isotoped so that $e_{b} \cup e_{+}$ is again level, but 
the width of $\Ggg$ has been reduced by at least $4m + 4$.
\end{lemma}

\begin{proof} The cycle $e_{+}$ divides $P$ into two disks $P_{1}
\cup P_{2}$.  Without loss of generality, assume that $e_{b}$ lies in
$P_{2}$.  Let $S_{P,i}^{Q}, i = 1, 2$ denote the component of
$S_{P}^{Q}$ lying above $P_{i}$.

\bigskip

{\bf Case 1:} Some maximum (resp. minimum) of $e_{-}$ can be pushed 
down (resp. up) past $P$.

Note that a plane just above or below $P$ intersects $e_{-}$ in at
least $2m$ points.  If the maximum that is pushed down is not the
maximum contiguous to the end of $e_{b}$ then the move instantly
reduces the width of $\Ggg$ by $8$, per \ref{exam:level}.  More
importantly, after the move $\Ggg$ is in a position to apply Lemma
\ref{lemma:level}, and so we can reduce the width by at least a
further $4(2m - 2)$.  Thus the total width is reduced by at least $8m
\geq 4m + 4$.  

If the maximum that is pushed down is contiguous to the end of 
$e_{b}$, the effect on width is to first push a regular maximum down 
past a $Y$-vertex (on $e_{+}$) and then to convert the regular maximum 
and the $Y$-vertex on $e_{-}$ into a single $\lll$ vertex on 
$e_{-}$.  The first move reduces the width by $4$ and the second move 
(eliminating a critical point) reduces it by at least a further 
$4m + 2$.

\bigskip

{\bf Case 2:} Some maximum of $e_{-}$ lies in
$S_{P,1}^{Q}$.  

The descending disk of any maximum in this region can't
intersect the end $C \cap Q$, since that end is too high.  Hence the
intersection of such a descending disk with $C$ consists entirely of
components that are inessential in $C$.  It follows that a disk in
$S_{P,1}^{Q}$ can be found that isotopes a maximum of $e_{-}$ in
$S_{P,1}^{Q}$ down to a level below $P$, returning us to Case 1.

\bigskip

Let $H$ be a regular neighborhood of $\Ggg$ and continue to call 
$P_{i}$ the disks obtained by removing the boundary collars given by 
$H \cap P_{i}$.  Then each $P_{i}$ is a disk punctured by meridians of 
$H$ associated with points on $e_{-}$.  Since $P$ was a dividing 
sphere for $e_{-}$, there are an odd number $p$ of such meridians (the 
point of $e_{-}$ at the end of $e_{b}$ does not, of course, give rise 
to such a meridian).  $P$ divides $\bdd H$ into $p+1$ components; 
$p-1$ of them are annuli $A_{1},\ldots,A_{p-1}$ lying between meridian 
disks associated to points in $e_{-} \cap P$.  Two components, 
$U_{\pm}$ are pairs of pants, with boundary of each consisting of 
$\bdd P_{1}, \bdd P_{2}$ and the boundary of a meridian associated to 
a point of $e_{-} \cap P$.  Choose notation so that $U_{+}$ lies above 
$P$, the meridian curves in $\bdd H$ associated to points of $e_{-} 
\cap P$ occur in order $\mu_{1},\ldots,\mu_{p}$ along $e_{-}$, with 
$\mu_{1} \subset \bdd U_{+}$ and $\mu_{p}\subset \bdd U_{-}$ and, 
finally, $\bdd A_{i} = \mu_{i} \cup \mu_{i+1}$.

Not surprisingly, we consider how a meridian disk $E$ of $S^{3} - H$
intersects $P$.  It will eventually be useful to have chosen $E$,
among all possible meridian disks, to minimize $|E \cap P|$.  Of course
if $E$ is disjoint from $P$ then its boundary can't be a meridian
curve of $e_{-}$ (every sphere in $S^{3}$ separates) so it must be
parallel to $\bdd P_{1}$.  But then it's easy to see that $\Ggg$ is in
fact planar, contradicting hypothesis.  If there are any closed
components of $E \cap (P_{1} \cup P_{2})$ then an innermost one on $E$
can be used to push a maximum below $P$ or a minimum above $P$.  Then
we are in Case 1 and the argument is complete.  A similar argument
applies if an outermost disk $E_{0}$ cut off from $E$ by $P_{i}$ is
incident to one of the $A_{i}$.  We conclude that $E \cap P$ consists
entirely of arcs and, furthermore, each outermost disk is incident
only to one of $U_{\pm}$.  Let $E_{0}$ be any such outermost disk, 
with boundary the union of two arcs $\aaa \subset \bdd\eta(H) \cap E$ and 
$\bbb \subset P \cap E$ in $E$. Consider the possibilities for $\aaa$.

\bigskip

{\bf Case 3:} One or both ends of $\aaa$ is incident to $\bdd P_{2}$.  

The other end of $\aaa$ can't be incident to $\bdd P_{1}$, for the arc
$\bbb = E_{0} \cap P$ lies either in $P_{1}$ or $P_{2}$.  If the
other end is incident to $\mu_{1}$ then it can be used to pull the
maximum of $e_{-}$ contiguous to the end of $e_{b}$ down below $P$,
again placing us in Case 1.  Similarly if the other end of $\aaa$ is
incident to $\mu_{p}$.  In fact, if both ends of $\aaa$ lie on $\bdd
P_{2}$ we can accomplish the same thing, essentially using $E_{0}$
much like a cap.  

\bigskip

{\bf Case 4:} Exactly one end of $\aaa$ is incident to $\bdd P_{1}$.  

Again, the other end of $\aaa$ can't be incident to $\bdd P_{2}$.  
Suppose it is incident to $\mu_{p}$.  Then, since an arc in a pair of 
pants is determined up to proper isotopy by its end points, the arc 
$\aaa$ runs once along the length of $e_{b}$, then over the minimum of 
$e_{-}$ that is adjacent to the end of $e_{b}$ and ends in $\mu_{p} 
\subset P_{1}$.  The disk $E_{0}$ can be used to slide $e_{b}$, 
keeping the end at $e_{+}$ fixed, until $e_{b}$ becomes the arc $\bbb 
\subset P_{1}$.  See Figure \ref{fig:flip}.  Afterwards, the width 
is unaffected, but all $m$ maxima now lie in the component 
$S_{P,2}^{Q}$ that no longer contains $e_{b}$.  In effect, we are in 
Case 2 and so we are finished once again.  The same argument applies 
if the other end of $\aaa$ is at $\mu_{1}$: Since the interior of 
$E_{0}$ is disjoint from $\Ggg$ the slide of $e_{b}$ to $\bbb$ has no 
effect on the maxima in $S_{P,2}^{Q}$, or on the cylinder $C$.  (The 
edge $e_{b}$ just passes through $C$).

\begin{figure}
\centering
\includegraphics[width=0.65\textwidth]{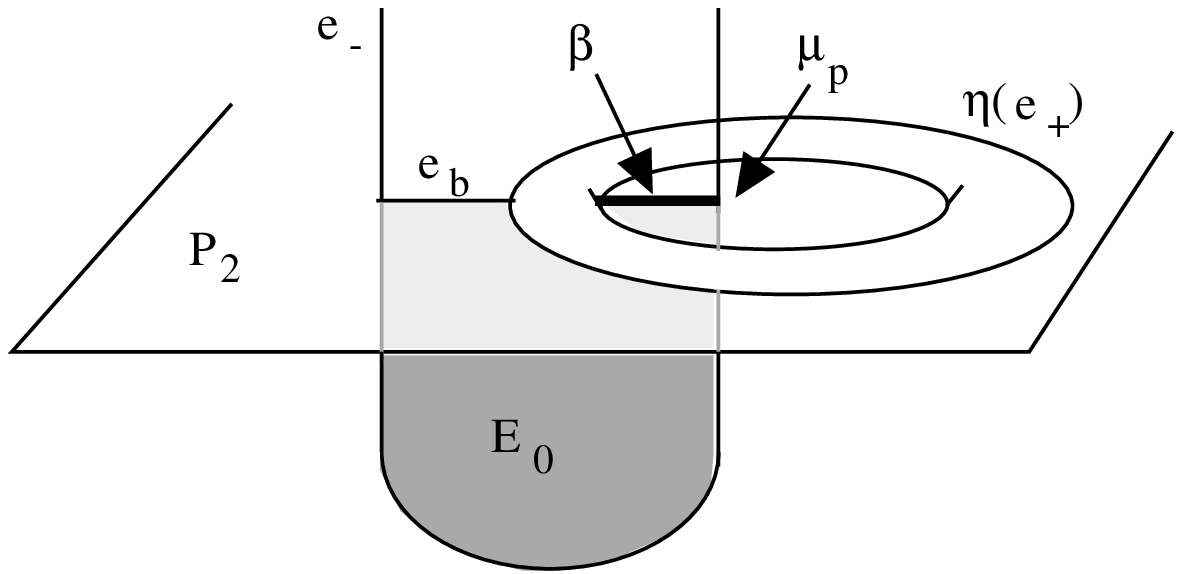}
\caption{Case 4}  \label{fig:flip}
\end{figure}

\bigskip

{\bf Case 5:} Both ends of $\aaa$ are incident to $\bdd P_{1}$.  

Suppose, to be concrete, that $E_{0}$ lies above $P$, so it forms a
kind of cap or shroud over the part $e_{0}$ of $e_{-}$ that lies
between $e_{b}$ and $\mu_{1}$.  Let $A$ denote the annulus half of
$\bdd \eta(e_{b} \cup e_{+})$ that lies above $P$ and let $P_{u}$
denote the plane $P_{1} \cup A \cup P_{2}$.  Then $\bdd E_{0} \subset
P_{u}$ consists of two arcs, $\aaa \subset A$ and $\bbb \subset
P_{1}$.  A descending disk for the maximum $e_{0}$ also has boundary
consisting of two arcs, one being $e_{0}$ itself and the other an arc
in $P_{u}$.  A standard innermost disk, outermost arc argument shows
that such a disk $D$ can be found disjoint from $E_{0}$, so $\bdd D$
lies in the disk in $P_{u}$ bounded by $\bdd E_{0}$.  In fact, $E_{0}$
can be used to remove (by piping to $E_{0}$ and then over it) any arc
of $\bdd D \cap A$ which is parallel to $\aaa$ in the punctured
annulus $A - e_{0}$.  The upshot is that, if we choose $D$ so that the
arc $\ddd = \bdd D \cap P_{u}$ intersects $A$ in a minimal number of
components, then in fact $\ddd$ consists of a single arc in $A$
(running from the end of $e_{0}$ to $\bdd P_{1}$) and a single arc in
$P_{1}$.  Once this is accomplished, the disk $D$ can be used instead
of $E_{0}$ in the proof of Case 4, completing the argument in this
case.

\bigskip

{\bf Case 6:} The general case.

Following cases 3 to 6, the only remaining case to consider is one in
which {\em every} outermost arc cut off by $P_{1} \cup P_{2}$ has both
ends incident to $\mu_{1}$ (when the disk it cuts off lies above $P$)
or both ends incident to $\mu_{p}$ (when the disk it cuts off lies
below $P$).  Notice that, in either case, the outermost arc forms a
loop in $P$ with both ends either at $\mu_{1}$ or at $\mu_{p}$.

\bigskip

{\bf Claim:}  For any $\mu_{i}, 1 \leq i \leq p$, there is an arc of 
$E \cap P$ forming a loop at $\mu_{i}$.  

The proof of the claim is a particularly easy application of outermost
forks.  Cf \cite{Sc} for details beyond this sketch: Label the ends of
arcs of $P \cap E$ in $\bdd E$ that lie on the meridians $\mu_{1},
\ldots, \mu_{p}$ by the number of the corresponding meridian.  We have
just demonstrated that each outermost arc has either both ends
labelled $1$ or both ends labelled $p$.  To the collection of arcs $E
\cap P$ there is naturally associated a tree in $E$, with a vertex in
each component of $E - P$ and an edge connecting any vertices
corresponding to adjacent components.  Consider an outermost fork of
this tree.  Two adjacent tines of this fork have ends labelled $(1,1)$
or $(p, p)$.  In order to get from one labelling to the other, the arc
of $\bdd E$ that lies between the ends of the two adjacent tines must
go sequentially through every label from $1$ to $p$ (perhaps more than
once).  Since each arc of $E \cap P$ it passes by is parallel to an
outermost arc, its labels must be the same.  The result is a
collection of arcs containing all labels $1, \ldots, p$ and having the
same label at each end.  (See Figure \ref{fig:fork}).  These arcs,
when considered in $P$, form loops at every meridian $\mu_{p}$.

\begin{figure}
\centering
\includegraphics[width=0.65\textwidth]{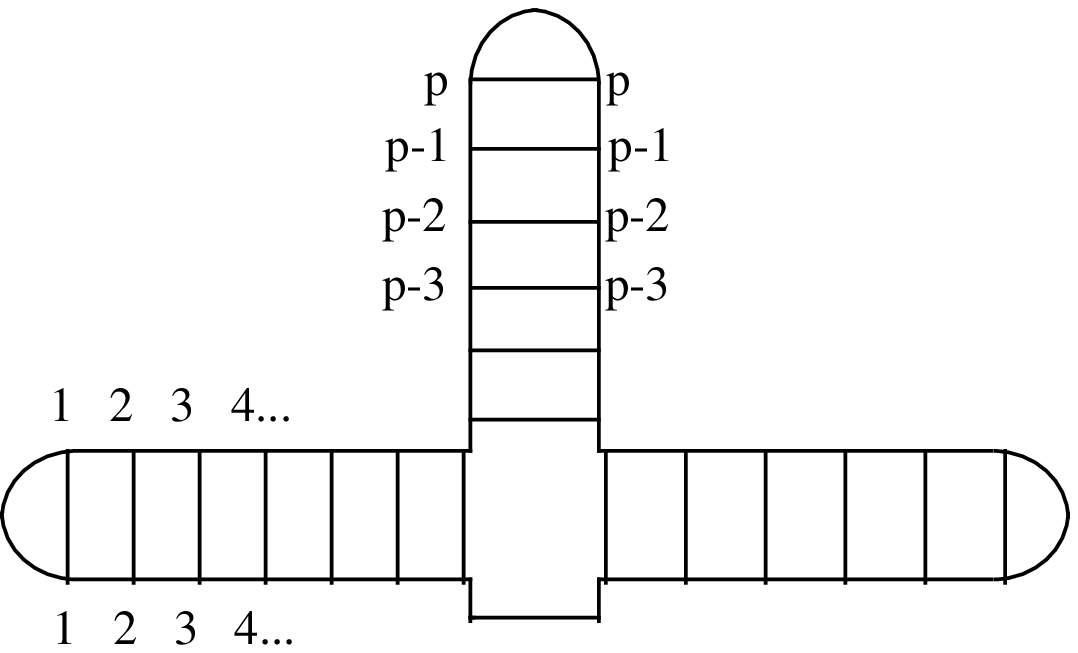}
\caption{}  \label{fig:fork}
\end{figure}

\bigskip

Having established the claim, consider this consequence: An innermost 
such loop contains no meridian in its interior.  This means that an 
innermost loop can be used to $\bdd$-compress $E$ to $\bdd H$, 
dividing $E$ into two disks, at least one of which is still a 
meridian 
disk and each of which intersects $P$ in fewer arcs.  Since $E$ was 
initially chosen to minimize $E \cap P$, this is impossible.
\end{proof}

\end{document}